\documentclass[3p]{elsarticle}
\usepackage{amsmath}
\usepackage{amsfonts,amssymb}

\newtheorem{theorem}{Theorem}{\bf}{\it}
{\bf}{\it}
{\bf}{\rm}
\newtheorem{lemma}{Lemma}{\bf}{\it}
\newtheorem{proposition}{Proposition}{\bf}{\it}
{\bf}{\it}

\newcommand{\zr}{\ltimes}
\newcommand{\Real}{\mathbb{R}}
\newcommand{\Co}{\mathbb{C}}

\newcommand{\g}{\mathfrak{g}}
\newcommand{\h}{\mathfrak{h}}

\newcommand{\so}{\mathfrak{so}}

\def\sl{\mathfrak{sl}}
\newcommand{\gl}{\mathfrak{gl}}

\def\u{\mathfrak{u}}
\def\sp{\mathfrak{sp}}
\newcommand{\f}{\mathfrak{f}}

\newcommand{\R}{\mathcal{R}}

\newcommand{\id}{\mathop\text{\rm id}\nolimits}
\newcommand{\spa}{\mathop\text{{\rm span}}\nolimits}
\newcommand{\Hom}{\mathop\text{\rm Hom}\nolimits}
\newcommand{\End}{\mathop\text{\rm End}\nolimits}

\newcommand{\ad}{\mathop\text{\rm ad}\nolimits}

\newcommand{\pr}{\mathop\text{\rm pr}\nolimits}
\newcommand{\Op}{\mathop\text{\rm Op}\nolimits}

\def\g{\mathfrak{g}}
\def\h{\mathfrak{h}}
\def\so{\mathfrak{so}}

\def\sl{\mathfrak{sl}}
\def\sp{\mathfrak{sp}}
\def\gl{\mathfrak{gl}}

\def\u{\mathfrak{u}}

\def\spa{\mathop\text{{\rm span}}\nolimits}
\def\Hom{\mathop\text{\rm Hom}\nolimits}
\def\pr{\mathop\text{\rm pr}\nolimits}

\def\End{\mathop\text{\rm End}\nolimits}
\def\id{\mathop\text{\rm id}\nolimits}

\def\ad{\mathop\text{\rm ad}\nolimits}

\def\Im{\mathop\text{\rm Im}\nolimits}

\def\E{{\cal E}}
\def\F{{\cal F}}

\def\zr{\ltimes}

\def\spa{\text{\rm span}}
\def\Hom{\text{\rm Hom}}

\def\End{\text{\rm End}}

\def\Real{\mathbb{R}}
\def\Co{\mathbb{C}}
\def\g{\mathfrak{g}}
\def\h{\mathfrak{h}}

\def\so{\mathfrak{so}}

\def\sl{\mathfrak{sl}}
\def\gl{\mathfrak{gl}}

\def\u{\mathfrak{u}}

\def\sp{\mathfrak{sp}}
\def\a{\mathfrak{a}}
\def\f{\mathfrak{f}}

\def\R{\mathcal{R}}

\def\B{\mathcal {B}}
\def\D{\mathcal {D}}

\def\K{\mathcal {K}}
\def\H{\mathbb {H}}

\def\F{\mathcal {F}}
\def\id{\text{\rm id}}

\def\pr{\text{\rm pr}}

\def\im{\text{\rm Im}}
\def\re{\text{\rm Re}}

\def\ad{\text{\rm ad}}

\def\Mat{\mathop\text{\rm Mat}\nolimits}

\begin{document}

\begin{frontmatter}

\title{Holonomy algebras of pseudo-quaternionic-K\"ahlerian manifolds \\ of signature $(4,4)$}

\author{Natalia I. Bezvitnaya}
\address{Department of Mathematics and Statistics, Faculty of Science, Masaryk University in Brno,
Kotl\'a\v rsk\'a~2, 611~37 Brno, Czech Republic\\
\ead{bezvitnaya@math.muni.cz}}

\begin{abstract} Possible holonomy algebras of pseudo-quaternionic-K\"ahlerian manifolds of signature
$(4,4)$ are classified. Using this, a new proof of the
classification of simply connected pseudo-quaternionic-K\"ahlerian
symmetric spaces of signature $(4,4)$ is obtained.
\end{abstract}


\end{frontmatter}

\section{Introduction}
The classification of  holonomy algebras of Riemannian manifolds
is well known and it has a lot of applications both in geometry
and physics, see e.g. \cite{Ber,Besse,Bryant,Jo,Gibbons09}. Lately
the theory of pseudo-Riemannian geometries has been steadily
developing. In particular, a classification of holonomy algebras
of pseudo-Riemannian manifolds is an actual problem of
differential geometry. It is solved only in some cases. Holonomy
algebras of 4-dimensional Lorentzian manifolds are classified in
\cite{Schell}. Classification of holonomy algebras of Lorentzian
manifolds  is  obtained in \cite{BB-I,Le,Gal5,ESI}; classification
of holonomy algebras
  of pseudo-K\"ahlerian manifolds of index 2 is achieved in  \cite{GalDB}.
 These algebras are contained in
 $\so(1,n+1)$ and $\u(1,n+1)\subset\so(2,2n+2)$, respectively.
In \cite{I99} a partial classification of holonomy algebras of
pseudo-Riemannian manifolds of signature $(2,n)$ is obtained. In
\cite{BI97} holonomy algebras contained in $\so(n,n)$ and living
invariant two $n$-dimensional  isotropic complementary subspaces
are considered. In \cite{BI97,GT01} holonomy algebras of
pseudo-Riemannian manifolds of signature $(2,2)$ are classified.
More details can be found in the recent review \cite{IRMA}.

 As the next step, in this paper we begin to study the holonomy algebras
 contained in $\sp(1)\oplus\sp(1,n+1)\subset\so(4,4n+4)$. Here we
 consider the case $n=0$. The technic developed in the present
 paper together with the results from \cite{B1} will allow to get
 the classification for arbitrary $n$, this will be done in
 another paper.

Recall that a pseudo-quaternionic-K\"ahlerian manifold  is a
pseudo-Riemannian manifold $(M,g)$ of signature $(4r,4s)$, $r+s>1$
together with a parallel quaternionic structure $Q\subset\so(TM)$,
i.e. three-dimensional linear Lie algebra $Q$ with a local basis
$I_{1},I_{2},I_{3}$ which satisfies the relations
$I_{1}^{2}=I_{2}^{2}=I_{3}^{2}=-\id$,
$I_{3}=I_{1}I_{2}=-I_{2}I_{1}$. The holonomy algebra $\g$ of such
manifold is contained in $\sp(1)\oplus\sp(r,s)$. Conversely, any
simply connected pseudo-Riemannian manifold with such holonomy
algebra is pseudo-quaternionic-K\"ahlerian.

Quaternionic-K\"ahlerian geometry is of increasing interest both
in mathematics and mathematical physics, see e.g.
\cite{AC01,AC05,KO07,March,Salamon,Swann}.

 A pseudo-hyper-K\"ahlerian manifold is a pseudo-quaternionic-K\"ahlerian manifold
 such that the endomorphisms $I_{1},I_{2},I_{3}$ are defined globally and parallel.
 In this case, the holonomy algebra
$\g$ is contained in  $\sp(r,s)$.

In \cite{AC05} it is proved that the curvature tensor $R$ of a
pseudo-quaternionic-K\"ahlerian manifold of signature $(4r,4s)$
can be written as
\begin{equation}\label{RR0W}R(X,Y)=\nu R_0+{\cal W},\end{equation}
where $\nu=\frac{{\rm scal}}{4m(m+2)}$ ($m=r+s$) is the reduced
scalar curvature,
\begin{equation}\label{R0}R_0(X,Y)=\frac{1}{2}\sum_{\alpha=1}^3g(X,I_\alpha
Y)I_\alpha-\frac{1}{4}\left(X\wedge Y+\sum_{\alpha=1}^3I_\alpha
X\wedge I_\alpha Y\right),\end{equation} $X,Y\in TM$, is the
curvature tensor of the quaternionic projective space
$\mathbb{PH}^{r,s}$, and ${\cal W}$ is an algebraic curvature
tensor with zero Ricci tensor. It is proved that  ${\rm scal}\neq
0$ if and only if the holonomy algebra contains $\sp(1)$. This
shows that a pseudo-quaternionic-K\"ahlerian manifold $(M,g)$ is
pseudo-hyper-K\"ahlerian if and only if its scalar curvature is
zero. In \cite{AC05} it is proved also that any
pseudo-quaternionic-K\"ahlerian manifold with non-zero scalar
curvature is locally indecomposable, i.e. it is not locally a
product of two pseudo-Riemannian manifolds of positive dimension.

We prove the following theorem\footnote{In the first version of
this paper the first possibility for $\g$ in Theorem \ref{Thind}
was missed, this is pointed out by Bastian Brandes. The mistake
was caused by the false statement that if $\h\subset\sp(4,\Real)$
is an irreducible subalgebra with a non-trivial first
prolongation, then $\h=\sp(4,\Real)$. In fact,
$\h=\sl(2,\Co)=\sp(2,\Co)\subset\sp(4,\Real)$ is the second
possibility for $\h$.}.

\begin{theorem}\label{Thind} Let $(M,g)$ be a non-flat pseudo-quaternionic-K\"ahlerian manifold
of signature $(4,4)$. Then $(M,g)$ is locally indecomposable. If
the holonomy algebra $\g\subset\sp(1)\oplus \sp(1,1)$ of $(M,g)$
is not irreducible then one of the following holds:
\begin{description}
\item[1)] there exists a basis $p_1,...,p_4,q_1,...,q_4$ of
$\Real^{4,4}$ such that the metric on $\Real^{4,4}$ has the Gram
matrix $\left(\begin{array}{cc}0&E_4\\E_4&0\end{array}\right)$, it
holds
\begin{equation}\label{poterg}
\g=\left.\left\{\Op\left(\begin{array}{cc}A&0\\0&-A^t\end{array}\right)\right|A\in\sl(2,\Co)
\subset\sp(4,\Real)\right\}\subset\sp(1,1),\end{equation} and the
hyper-complex structure on $\Real^{4,4}$ can be chosen in the
following way:
$$I_1=\left(\begin{array}{cc}I&0\\0&-I\end{array}\right), \quad
I_2=\left(\begin{array}{cc}0&J\\J&0\end{array}\right),\quad
I_3=I_1I_2=\left(\begin{array}{cc}0&K\\K&0\end{array}\right),$$
where $$
I=\left(\begin{array}{cc}0&-E_2\\E_2&0\end{array}\right),\quad
J=\left(\begin{array}{cccc}0&1&0&0\\-1&0&0&0\\0&0&0&-1\\0&0&1&0
\end{array}\right),\quad K=\left(\begin{array}{cccc}0&0&0&-1\\0&0&1&0\\0&-1&0&0\\1&0&0&0
\end{array}\right);$$
\item[2)]  $\g$  preserves a $4$-dimensional isotropic
$I_1,I_2,I_3$-invariant subspace $W\subset\Real^{4,4}$.
\end{description}\end{theorem}

Here the pseudo-Euclidean space $\Real^{4,4}$ is considered as the
tangent space at some point of the manifold~$(M,g)$.

Suppose that $\g$ preserves a 4-dimensional
$I_{1},I_{2},I_{3}$-invariant isotropic subspace
$W\subset\mathbb{R}^{4,4}$. We identify $\mathbb{R}^{4,4}$ with
the pseudo-quaternionic Hermitian space $\mathbb{H}^{1,1}$. Then
$W\subset\mathbb{H}^{1,1}$ is an isotropic quaternionic line. We
fix  isotropic vectors $p\in W$ and $q\in\mathbb{H}^{1,1}$  such
that $g(p,q)=1$, then $W=\mathbb{H}p$ and
$\mathbb{H}^{1,1}=\mathbb{H}p\oplus\mathbb{H}q$. Denote by
$\sp(1,1)_{\mathbb{H}p}$ the maximal subalgebra of $\sp(1,1)$ that
preserves the quaternionic isotropic line $\mathbb{H}p$, this Lie
algebra has the matrix form:
\begin{equation}\label{spHp}\sp(1,1)_{\H p}=
\left.\left\{\Op\left(\begin{array}{cc}a&b\\0&-\bar
a\end{array}\right)\right|a\in\H, \,b\in\Im
\H\right\}.\end{equation} Here $\Op A$ denotes the $\H$-linear
endomorphism of $\H^{1,1}$ given by a matrix $A$, see Section
\ref{secPril}.

We prove the following two theorems.

\begin{theorem}\label{Th2} Let $(M,g)$ be a non-flat pseudo-hyper-K\"ahlerian manifold
of signature $(4,4)$. Suppose that the holonomy algebra
$\g\subset\sp(1,1)$ of $(M,g)$ preserves a $4$-dimensional
isotropic $I_1,I_2,I_3$-invariant subspace $W\subset\Real^{4,4}$.
Then $\g$ is conjugated to one of the following  subalgebras of
$\sp(1,1)_{\H p}$:
\begin{align*}\g_1&=\sp(1,1)_{\H p},\qquad\qquad\qquad
\g_2=\left.\left\{\Op\left(\begin{array}{cc}a&b\\0&-\bar a\end{array}\right)\right|a,b\in\Im \H\right\},\\
\g_3&=\left.\left\{\Op\left(\begin{array}{cc}a&b\\0&-\bar
a\end{array}\right)\right|a\in\Real\oplus\Real i,\,b\in\Im
\H\right\},\qquad
\g_4=\left.\left\{\Op\left(\begin{array}{cc}0&b\\0&0\end{array}\right)\right|b\in\Im \H\right\},\\
\g_5&=\left.\left\{\Op\left(\begin{array}{cc}a&b\\0&-\bar
a\end{array}\right)\right|a\in\Real\oplus\Real k,\,b\in\Real
i\oplus\Real j\right\},\quad
\g_6=\left.\left\{\Op\left(\begin{array}{cc}0&b\\0&0\end{array}\right)\right|b\in\Real
i\oplus\Real j\right\}.
\end{align*}
\end{theorem}

\begin{theorem}\label{Th3} Let $(M,g)$ be a  pseudo-quaternionic-K\"ahlerian manifold
of signature (4,4) with non-zero scalar curvature. Suppose that
the holonomy algebra $\g\subset\sp(1)\oplus\sp(1,1)$ of $(M,g)$
preserves a 4-dimensional isotropic $I_1,I_2,I_3$-invariant
subspace $W\subset\Real^{4,4}$. Then $\g$ is conjugated to one of
the following  subalgebras of $\sp(1)\oplus\sp(1,1)_{\H p}$: $$
\sp(1)\oplus\sp(1,1)_{\H p},\qquad
\h_0=\sp(1)\oplus\left.\left\{\Op\left(\begin{array}{cc}a&0\\0&-\bar
a\end{array}\right)\right|a\in\H\right\}.$$
\end{theorem}
The last Lie algebra can be written as
$\h_0=\left.\left\{\left(\begin{array}{cc}A&0\\0&-A^t\end{array}\right)\right|A\in
\mathfrak{co}(4)\right\}.$

To prove these theorems we use the fact that a holonomy algebra
$\g\subset\sp(1)\oplus\sp(1,1)$ is {a Berger algebra}, i.e. $\g$
is spanned by the images of the algebraic curvature tensors
$R\in\R(\g)$ of type $\g$. Recall that $\R(\g)$ is the space of
linear maps from $\wedge^2\Real^{4,4}$ to $\g$ satisfying the
first Bianchi identity. We find all subalgebras
$\g\subset\sp(1,1)$, find the spaces $\R(\g)$ and
$\R(\sp(1)\oplus\g)$ and check which of these algebras are Berger
algebras. Remark that in the above two theorems only possible
holonomy algebras are listed. We do not know if all these algebras
may appear as the holonomy algebras, to show this examples of
manifolds must be constructed.

In Section \ref{secR} we give the spaces $\R(\g)$ for each algebra
$\g$ from Theorems \ref{Th2} and \ref{Th3}.

In \cite{AC05} a classification of simply connected
pseudo-quaternionic-K\"ahlerian symmetric spaces of non-zero
scalar curvature is obtained. In \cite{AC01,KO07,KO08} simply
connected pseudo-hyper-K\"ahlerian symmetric spaces of index 4 are
classified.

In Section \ref{secSym} we show that if the manifold $(M,g)$ is
locally symmetric, then its holonomy algebra is one of the
following:
$$\left.\left\{\Op\left(\begin{array}{cc}0&b\\0&0\end{array}\right)\right|b\in\Im
\H\right\},\quad
\left.\left\{\Op\left(\begin{array}{cc}0&b\\0&0\end{array}\right)\right|b\in\Real
i\oplus\Real j\right\},\quad
\sp(1)\oplus\left.\left\{\Op\left(\begin{array}{cc}a&0\\0&-\bar
a\end{array}\right)\right|a\in\H\right\}.$$ Using this, we get a
new proof of the classification of pseudo-quaternionic-K\"ahlerian
symmetric spaces of signature (4,4) with non-irreducible holonomy
algebras. We give explicitly the curvature tensors of the obtained
spaces.

{\it Acknowledgement.} I am grateful to D.~V.~Alekseevsky and Jan
Slov\'ak for useful discussions, help and support. I am thankful
to Bastian Brandes for pointing out a mistake in the first version
of this paper. The author has been supported by the grant GACR
201/09/H012.

\section{Preliminaries}\label{secPril}
Let $\mathbb{H}^{m}$ be an $m$-dimensional quaternionic vector
space. A pseudo-quaternionic-Hermitian metric $g$ on
$\mathbb{H}^{m}$ is a non-degenerate $\mathbb{R}$-bilinear map
$g:\mathbb{H}^{m}\times\mathbb{H}^{m}\to\mathbb{H}$ such that
$g(aX,Y)=ag(X,Y)$ and $\overline{g(Y,X)}=g(X,Y)$, where
$a\in\mathbb{H}$, $X,Y\in\mathbb{H}^{m}$. Hence,
$g(X,aY)=g(X,Y)\bar{a}$. There exists a basis $e_{1},...,e_{m}$ of
$\mathbb{H}^{m}$ and integers $(r,s)$ with $r+s=m$ such that
$g(e_{t},e_{l})=0$ if $t \neq l$, $g(e_{t},e_{t})=-1$ if $1 \leq t
\leq p$ and $g(e_{t},e_{t})=1$ if $p+1 \leq t \leq m$. The pair
$(r,s)$ is called  the signature of $g$. In this situation we
denote $\mathbb{H}^{m}$ by $\mathbb{H}^{r,s}$. The realification
of $\mathbb{H}^{m}$ gives us the vector space $\mathbb{R}^{4m}$
with the quaternionic structure $(i,j,k)$. Conversely, a
quaternionic structure on $\mathbb{R}^{4m}$, i.e. a triple
$(I_{1},I_{2},I_{3})$ of endomorphisms of $\mathbb{R}^{4m}$ such
that $I_{1}^{2}=I_{2}^{2}=I_{3}^{2}=-\id$ and
$I_{3}=I_{1}I_{2}=-I_{2}I_{1}$, allows us to consider
$\mathbb{R}^{4m}$ as $\mathbb{H}^{m}$. A
pseudo-quaternionic-Hermitian metric $g$ on $\mathbb{H}^{m}$ of
signature $(r,s)$ defines on $\mathbb{R}^{4m}$ the
$i,j,k$-invariant pseudo-Euclidean metric $\eta$ of signature
$(4r,4s)$, $\eta(X,Y)=\re g(X,Y)$, $X,Y\in\mathbb{R}^{4m}$.
Conversely, a $I_{1},I_{2},I_{3}$-invariant pseudo-Euclidean
metric on $\mathbb{R}^{4m}$ defines a
pseudo-quaternionic-Hermitian metric $g$ on $\mathbb{H}^{m}$,
$$g(X,Y)=\eta(X,Y)+i\eta(X,I_{1}Y)+j\eta(X,I_{2}Y)+k\eta(X,I_{3}Y).$$

We will identify $1,i,j,k$ with $I_0,I_1,I_2,I_3$, respectively.
The identification $\Real^{4r,4s}\simeq\H^{r,s}$ allows to
multiply the vectors of $\Real^{4r,4s}$ by quaternionic numbers.

The Lie algebra $\sp(r,s)$ is defined as follows
\begin{align*} \sp(r,s)&=\{f\in\so(4r,4s) |
[f,I_{1}]=[f,I_{2}]=[f,I_{3}]=0\}\\&=
\{f\in\End(\mathbb{H}^{r,s})|\, g(fX,Y)+g(X,fY)=0 \text{ for all }
X,Y\in\mathbb{H}^{r,s}\}.\end{align*} Denote by $\sp(1)$ the
subalgebra in $\so(4r,4s)$  generated by the $\mathbb{R}$-linear
maps $I_{1},I_{2},I_{3}$.

Clearly, the tangent space of a pseudo-quaternionic-K\"ahlerian
manifold $(M,g)$ at a point $x\in M$ one can identify with
$(\mathbb{R}^{4r,4s},\eta_{x},(I_{1})_{x},(I_{2})_{x},(I_{3})_{x})=(\mathbb{H}^{r,s},g_{x})$.
Then the holonomy algebra of a  pseudo-quaternionic-K\"ahlerian
manifold is identified with a subalgebra
$\g\subset\sp(1)\oplus\sp(r,s)$.

Let $(V,\eta)$ be a pseudo-Euclidean space and $\g\subset\so(V)$
be a subalgebra. The space of curvature tensors $\mathcal{R}(\g)$
of type $\g$ is defined as follows
$$\mathcal{R}(\g)=\{R\in\Hom(\wedge^{2} V,\g)\ |\
R(u\wedge v)w+R(v\wedge w)u+R(w\wedge u)v=0 \ for \ all \
 u,v,w\in V\}.$$
 It is known that any
 $R\in\mathcal{R}(\g)$ satisfies \begin{equation}\label{(*)}
 \eta(R(u\wedge v)z,w)=\eta(R(z\wedge w)u,v)
 \end{equation} for all $u,v,w,z\in V$.

 Denote by $L(\mathcal{R}(\g))$ the vector subspace of $\g$
 spanned by the elements $R(u\wedge v)$ for all
 $R\in\mathcal{R}(\g)$ and $u,v\in V$.
{\it A subalgebra $\g\subset\so(r,s)$ is called a Berger algebra
if $L(\mathcal{R}(\g))=\g$.} From the Ambrose-Singer theorem it
follows that if $\g\subset\so(V)$ is the holonomy algebra of a
pseudo-Riemannian manifold, then $\g$ is a Berger algebra.
Therefore,  Berger algebras may be considered as the candidates to
the holonomy algebras.

Now we summarize some facts about quaternionic vector spaces. Let
$\mathbb{H}^{m}$ be an m-dimensional quaternionic vector space and
$e_{1},...,e_{m}$ a basis of $\mathbb{H}^{m}$. We identify an
element $X\in\mathbb{H}^{m}$ with the column $(X_{t})$ of the left
coordinates of $X$ with respect to this basis,
$X=\sum_{t=1}^{m}X_{t}e_{t}$. Let
$f:\mathbb{H}^{m}\to\mathbb{H}^{m}$ be an $\mathbb{H}$-linear map.
Define the matrix $\Mat_{f}$ of $f$ by the relation
$fe_{l}=\sum_{t=1}^{m}(\Mat_{f})_{tl}e_{t}$. Now if
$X\in\mathbb{H}^{m}$, then $fX=(X^{t}\Mat_{f}^{t})^{t}$ and
because of the non-commutativity of the quaternionic numbers this
is not the same as $\Mat_{f}X$. Conversely, to an $m\times m$
matrix $A$ of the quaternionic numbers we put in correspondence
the linear map $\Op{A}:\mathbb{H}^{m}\to\mathbb{H}^{m}$ such that
$\Op{A}\cdot X=(X^{t}A^{t})^{t}$. If
$f,g:\mathbb{H}^{m}\to\mathbb{H}^{m}$ are two $\mathbb{H}$-linear
maps, then $\Mat_{fg}=(\Mat_{g}^{t}\Mat_{f}^{t})^{t}$. Note that
the multiplications by the imaginary quaternionic numbers are not
$\mathbb{H}$-linear maps. Also, for $a,b\in\mathbb{H}$ holds
$\overline{ab}=\bar{b}\bar{a}$. Consequently for two square
quaternionic matrices we have
$(\overline{AB})^{t}=\bar{B}^t\bar{A}^t$.

Let $R\in\R(\sp(r,s))$. Using \eqref{(*)} it is easy to show that
for any  $1\leq \alpha\leq 3$ and $X,Y\in\Real^{4r,4s}$,
\begin{equation}R(I_\alpha X,Y)=-R(X,I_\alpha Y)\end{equation} holds.
Hence,
\begin{equation}R(xX,Y)=R(X,\bar xY)\end{equation}
for all $x\in \H$ and $X,Y\in\Real^{4r,4s}$.

Let $W\subset\mathbb{R}^{4,4}$ be an $I_{1},I_{2},I_{3}$-invariant
isotropic subspace. Then  $W$ may be seen as an isotropic line in
$\mathbb{H}^{1,1}$. Fix a nonzero vector $p\in W$ then
$W=\mathbb{H}p$ and $g(p,p)=0$. Let $q\in\mathbb{H}^{1,1}$ be any
vector such that $g(q,q)=0$ and $g(p,q)=1$. Obviously, such vector
exists. Any other $q'\in\mathbb{H}^{1,1}$ with this property has
the form $q'=\alpha p+q$, where $\alpha\in\Im\H$. We get the basis
$(p,q)$ of $\H^{1,1}$.

Denote by $\sp(1,1)_{\mathbb{H}p}$ the maximal subalgebra of
$\sp(1,1)$ that preserves the quaternionic isotropic line
$\mathbb{H}p$. This algebra has the matrix form \eqref{spHp}. One
can easily find that the Lie brackets are the following:
$$\left[\Op\left(\begin{array}{cc}a&0\\0&-\bar
a\end{array}\right),\Op\left(\begin{array}{cc}0&b\\0&0\end{array}\right)\right]=
\left(\begin{array}{cc}0&2\Im ba\\0&0\end{array}\right) ,$$
$$\left[\Op\left(\begin{array}{cc}a_1&0\\0&-\bar
a_1\end{array}\right),\Op\left(\begin{array}{cc}a_2&0\\0&-\bar
a_2\end{array}\right)\right]= \left(\begin{array}{cc}
a_{2}a_{1}-a_1a_2&0\\0&-(\overline{a_{2}a_{1}-a_1a_2})\end{array}\right)
,$$ where $a_1,a_2\in\H$ and $b\in\im\H$. We get the
decomposition:
$$\sp(1,1)_{\H p}=(\mathbb{R}\oplus\sl(1,\H))\ltimes\Im\mathbb{H}.$$
Note that $\sl(1,\H)=\Im\H$. We use this notation in order to
differ it from the last term in the above decomposition. Below for
the matrix $\left(\begin{array}{cc}a&b\\0&-\bar
a\end{array}\right)$ is using the notation $(a,b)$.

Coming back to $\mathbb{R}^{4,4}$ we get the basis
$p,I_{1}p,I_{2}p,I_{3}p,q,I_{1}q,I_{2}q,I_{3}q$ with respect to
this basis the Gram matrix of $\eta$ has the form
$\left(\begin{array}{cc}0& E_4\\ E_4&0\end{array} \right)$ and
$$\sp(1,1)_{\mathbb{H}p}=\left.\left\{\left(\begin{array}{cc}
\begin{array}{cccc} a_0&-a_1&-a_2&-a_3\\
a_1&a_0&a_3&-a_2\\
a_2&-a_3&a_0&a_1\\
a_3&a_2&-a_1&a_0\end{array}&
\begin{array}{cccc}
 0&-b_1&-b_2&-b_3\\
b_1&0&b_3&-b_2\\
b_2&-b_3&0&b_1\\
b_3&b_2&-b_1&0
\end{array}\\
0&
\begin{array}{cccc}
-a_0&-a_1&-a_2&-a_3\\
a_1&-a_0&a_3&-a_2\\
a_2&-a_3&-a_0&a_1\\
a_3&a_2&-a_1&-a_0\end{array} \end{array}\right)\right|\
a_{0},a_1,a_2,a_{3}, \ b_{1},b_{2},b_{3}\in\mathbb{R}\right\}.$$

\section{Proof of Theorem \ref{Thind}}

 Suppose that
$\g\not\subset\sp(1,1)$. In \cite{AC05} it is proved that $(M,g)$
is locally indecomposable. If $\g$ is not irreducible, then $\g$
preserves a proper degenerate subspace  $W\subset\Real^{4,4}$.
Obviously, $W$ is $I_1,I_2,I_3$-invariant. Consequently, $\dim
W=4$. Note that $W\cap W^\bot$ is non-trivial, isotropic,
$I_1,I_2,I_3$-invariant and preserved by $\g$, hence $W\cap
W^\bot$ has dimension $4$, i.e. $W=W\cap W^\bot$ and $W$ is
isotropic.

Suppose now that $\g\subset\sp(1,1)$.

\begin{lemma}\label{Lemind1} If $\g\subset\sp(1,1)$ is
the holonomy algebra and $\g$ preserves a non-degenerate 4-dimensional vector subspace $W\subset\Real^{4,4}$
such that $W\cap I_s W=0$ for some $s\in\{1,2,3\}$, then $\g=0$.\end{lemma}

{\bf Proof.} Since $W$ is non-degenerate, we get the
$\g$-invariant orthogonal decomposition $\Real^{4,4}=W\oplus
W^\bot$. Then $\g=\h_1\oplus\h_2$, where $\h_1\subset\so(W)$ and
$\h_2\subset\so(W^\bot)$ are pseudo-Riemannian holonomy algebras.
For any $x\in W$, $I_s(x)$ can be uniquely written as
$I_s(x)=\varphi_1(x)+\varphi_2(x)$, where $\varphi_1:W\to W$ and
$\varphi_2:W\to W^\bot$ are linear maps. Since $W\cap I_s W=0$,
$\varphi_2$ is an injective, it is clear that it is an
isomorphism. Let $\xi\in \h_1$.
 The condition $\xi I_s=I_s\xi$ implies $\xi\varphi_1(x)=\varphi_1(\xi x)+\varphi_2(\xi x)$.
  Hence, $\varphi_2(\xi x)=0$. Since $\varphi_2$ is an injective, $\xi x=0$ for all $\xi\in\h_1$ and $x\in W$.
 Consequently, $\h_1=0$. Let $\xi\in \h_2$, $x\in W$.  We have $\xi I_s(x)=I_s\xi(x)=0$, therefore $\xi\varphi_2(x)=0$.
 Since $\varphi_2$ is an isomorphism, $\xi W^\bot=0$. Thus, $\xi=0$ and then $\h_2=0$. \qed

Suppose that $\g$ preserves a non-degenerate vector subspace
$W\subset\Real^{4,4}$. Let $W_1=W\cap I_1 W$. If $\dim W<4$, then
taking $W^\bot$ instead of $W$, we get that $W_1\neq 0$. If $\dim
W=4$, then from Lemma \ref{Lemind1} it follows that $W_1\neq 0$.
Let $W_2=W_1\cap I_2 W_1$. By the same arguments, we have that
$W_2\neq 0$. Hence $W_2\subset\Real^{4,4}$ is a non-degenerate
$I_1,I_2,I_3$-invariant subspace and $\dim W_2=4$. We get that
$\g\subset\sp(W_2)\oplus\sp(W_2^\bot)=\sp(1)\oplus\sp(1)$. Since
$\sp(1)\subset\so(4)$ is not a Riemannian holonomy algebra,
$\g=0$.

Thus, if $\g\subset\sp(1)\oplus\sp(1,1)$ is a non-trivial holonomy
algebra, then it does not preserve any non-degenerate subspace of
$\Real^{4,4}$ (such algebra is called weakly-irreducible).

\begin{lemma}\label{Lemind2} If $\g\subset\sp(1,1)$ is a holonomy algebra and $\g$
preserves an isotropic 4-dimensional vector subspace
$W\subset\Real^{4,4}$ such that $W\cap I_s W=0$ for some
$s\in\{1,2,3\}$, then either $\g=0$ or $\g$ is the algebra 1) from
the statement of the theorem.
\end{lemma}

{\bf Proof.} Without loss of generality we may assume that $W\cap
I_2 W=0$. Let $p_1,...,p_4$ and $q_1,...,q_4$ be bases of $W$ and
$I_2W$, respectively, such that $\eta(p_i,q_j)=\delta_{ij}$. Then
the metric $\eta$ has the  Gram matrix
$\left(\begin{array}{cc}0&E_4\\E_4&0\end{array}\right)$.
Consequently,
\begin{equation}\label{gsub}\g\subset\left.\left\{
\left(\begin{array}{cc}A&0\\0&-A^t\end{array}\right)\right|A\in\gl(4,\Real)\right\}.\end{equation}
Since $I_2$ is $\eta$-invariant, we get $I_2=
\left(\begin{array}{cc}0&C\\B&0\end{array}\right),$ where $B$ and
$C$ are  $4\times 4$-squire matrices such that $B^t=-B$ and
$C^t=-C$. Since $I_2^2=-\id$, we get that $CB=BC=-E_4$, hence,
$C=-B^{-1}$. Since $\g$ commutes with $I_2$, for any matrix from
\eqref{gsub} holds $A^tB+BA=0$. Let $\omega$ be the skew-symmetric
bilinear form on $W$ given by the matrix $B$ with respect to the
basis $p_1,...,p_4$. Then, $\g|_{W}\subset\sp(W,\omega)$. If
$\g\neq 0$, then since $\g\subset\sp(1,1)$ is weakly-irreducible,
$\g|_{W}\subset\sp(W,\omega)$ must be irreducible. Let
$R\in\R(\g)$. It is easy to see that $R(X,Y)=0$ whenever $X,Y\in
W$ or $X,Y\in I_2 W$. Next, for each fixed $Y\in I_2 W$,
$R(Y,\cdot|_{W})|_W$ belongs to the first prolongation
$(\g|_{W})^{(1)}$ of the subalgebra
$\g|_{W}\subset\sp(W,\omega)=\sp(4,\Real)$. If
$(\g|_{W})^{(1)}=0$, then  $\R(\g)=0$, i.e. $\g=0$. The only
irreducible subalgebras of $\sp(4,\Real)$ with non-trivial first
prolongations are $\sp(4,\Real)$ and
$\sp(2,\Co)\subset\sp(4,\Real)$ \cite{Bryant}. Thus $\g|_{W}$ is
either one of $\sp(4,\Real)$, $\sp(2,\Co)$, or $\g=0$.

Suppose that $\g|_{W}=\sp(2,\Co)$. Let $I$ denote the
corresponding complex structure on $W$. The vectors $p_1,...,p_4$
may be chosen in such a way that $p_3=Ip_1$, $p_4=Ip_2$ and $B=J$,
where $J$ is as in the statement of the theorem. Choose the
vectors $q_1,...,q_4$ such that $q_3=Iq_1$, $q_4=Iq_2$ and
$\eta(p_i,q_j)=\delta_{ij}$. Recall that the hyper-complex
structure $(I_1,I_2,I_3)$ must commute with $\g$. Note that if the
elements $I_1,I_2,I_3$ are not fixed, then they are not defined
canonically, while only $\spa_\Real\{I_1,I_2,I_3\}$ is defined
canonically.
 Let $\left(\begin{array}{cc}D&H\\F&-D^t\end{array}\right)\in\so(4,4)$,
then $F^t=-F$ and $H^t=-H$. This element commutes with $\g$ if and
only if $DA=AD$, $A^tF+FA=0$ and $AH+HA^t=0$, for all
$A\in\sp(2,\Co)\subset\sp(4,\Real)$. Clearly, $D$ is a linear
combination of $E_4$ and $I$, and $F$ and $H$ are linear
combinations of $J$ and $K$. We already know that
$I_2=\left(\begin{array}{cc}0&J\\J&0\end{array}\right)$. There
exist $a_1,...,a_6\in\Real$ such that
$I_1=\left(\begin{array}{cc}a_1E_4+a_2I&a_5J+a_6
K\\a_3J+a_4K&-a_1E_4+a_2I\end{array}\right)$. The condition
$I_1I_2=-I_2I_1$ implies $a_3=a_5=0$ and $a_4=a_6$. The condition
$I_1^2=-\id$ implies $a_1=0$, $a_2a_4=0$ and $a_1^2+a_4^2=1$.
Taking $a_2=1$ and $a_4=0$, we get that $I_1,I_2,I_3$ are as in
the statement of the theorem. Obviously, the other choose of $a_2$
and $a_4$ will define the same $\spa_\Real\{I_1,I_2,I_3\}$.

Further, we may assume that either $\g|_W=\sp(4,\Real)$ or $\g=0$.

Consider now $I_1$. Suppose that $W\cap I_1 W\neq 0$. If $\dim
W\cap I_1 W=2$, then $\g|_{W}$ preserves $W\cap I_1 W$,
consequently $\g=0$. If $\dim W\cap I_1 W=4$, i.e. $I_1$ preserves
$W$, then $\g|_{W}$ commutes with $I_1|_{W}$, hence  $\g=0$.

Suppose that $W\cap I_1 W=0$. Then there exist two isomorphisms
$\varphi_1$ and $\varphi_2$ of $W$ such that
$I_1(x)=\varphi_1(x)+I_2\varphi_2(x)$ for all $x\in W$. Since $\g$
commutes with $I_1$, $\g|_W$ commutes with both $\varphi_1$ and
$\varphi_2$. Suppose that $\g\neq 0$, then $\g|_{W}=\sp(W,\omega)$
and $\varphi_1$ and $\varphi_2$ are multiples of the identity. The
condition $I_1^2=-\id$ implies
$$-x=I_1^2(x)=I_1(\varphi_1(x)+I_2\varphi_2(x))=I_1\varphi_1(x)+I_1I_2\varphi_2(x)=
\varphi_1^2(x)+\varphi_2^2(x)+I_2[\varphi_2,\varphi_1](x)=\varphi_1^2(x)+\varphi_2^2(x).$$
This gives a contradiction. We get that $\g=0$. \qed

Suppose that $\g$ preserves a degenerate vector subspace
$W\subset\Real^{4,4}$.
 Let $W_1=W\cap I_2 W$.
If $\dim W<4$, then taking $W^\bot$ instead of $W$, we get that
$W_1\neq 0$. If $\dim W=4$, then from Lemma \ref{Lemind2} it
follows that either $W_1\neq 0$ or $\g$ satisfies condition 1) of
the statement of the theorem. Suppose that $W_1\neq 0$. Let
$W_2=W_1\cap I_1 W_1$. By the same arguments, we may assume that
$W_2\neq 0$. Note that by the above arguments, all proper
$\g$-invariant subspaces are degenerate. Hence
$W_2\subset\Real^{4,4}$ is a degenerate $I_1,I_2,I_3$-invariant
subspace and $\dim W_2=4$. Since $W_2\cap W_2^\bot\neq 0$ and it
is a $I_1,I_2,I_3$-invariant subspace, it coincides with $W_2$,
i.e. $W_2$ is isotropic. The theorem is proved. \qed

\section{Proof of Theorem \ref{Th2}} Let us find the space of curvature tensors
$\mathcal{R}(\sp(1,1)_{\H p})$ for the Lie algebra $\sp(1,1)_{\H
p}$.

Using the form $\eta$, the Lie algebra $\so(4,4)$ can be
identified with the space
$$\wedge^{2}\mathbb{R}^{4,4}=\spa\{u\wedge v=u\otimes v-v\otimes u
|u,v\in\mathbb{R}^{4,4}\}$$ in such a way that $(u\wedge
v)w=\eta(u,w)v-\eta(v,w)u$ for all $u,v,w\in \mathbb{R}^{4,4}$.
The metric $\eta$ defines the metric $\eta\wedge\eta$ on
$\wedge^{2}\mathbb{R}^{4,4}$. One can check that the element
$\Op\left(\begin{array}{cc}a&b\\0&-\bar
a\end{array}\right)\in\sp(1,1)_{\H p}$ corresponds to the bivector
\begin{align*}
&-a_0(p\wedge q+ip\wedge iq+ jp\wedge jq+ kp\wedge kq)+a_1(p\wedge
iq-ip\wedge q+ kp\wedge jq- jp\wedge kq)\\&+ a_2(-jp\wedge
q+p\wedge jq- kp\wedge iq+ ip\wedge kq)+a_3(-kp\wedge q+jp\wedge
iq-ip\wedge jq+ p\wedge kq)\\&+ b_1(p\wedge ip-jp\wedge
kp)+b_2(p\wedge jp+ip\wedge kp)+b_3(p\wedge kp-ip\wedge
jp).\end{align*}

Thus we get the decomposition
\begin{align*}
&\sp(1,1)_{\H p}=\Big(\mathbb{R}(p\wedge q+ip\wedge iq+ jp\wedge
jq+ kp\wedge kq)\oplus \mathbb{R}(p\wedge iq-ip\wedge q+ kp\wedge
jq- jp\wedge kq)\\
&\oplus \mathbb{R}(-jp\wedge q+p\wedge jq- kp\wedge iq+ ip\wedge
kq)\oplus\mathbb{R}(-kp\wedge q+jp\wedge iq-ip\wedge jq+ p\wedge
kq)\Big)\\&\zr\Big(\mathbb{R}(p\wedge ip-jp\wedge
kp)+\mathbb{R}(p\wedge jp+ip\wedge kp)+\mathbb{R}(p\wedge
kp-ip\wedge jp)\Big).
\end{align*}

Let $\g\subset\sp(1,1)_{\H p}$ be a subalgebra,
$R\in\mathcal{R}(\g)$. Using the above identification, $R$ can be
considered as the map
$$R:\wedge^{2}\mathbb{R}^{4,4}\to\g\subset\so(4,4)\simeq\wedge^2\mathbb{R}^{4,4}.$$
From \eqref{(*)}, we obtain
\begin{equation}\label{symR}\eta\wedge\eta(R(u\wedge v), z\wedge
w)=\eta\wedge\eta(R(z\wedge w), u\wedge v)\end{equation} for all
$u,v,z,w\in\mathbb{R}^{4,4}$. This shows that $R$ is a symmetric
linear map. Consequently $R$ is zero on the orthogonal complement
to $\g$ in $\wedge^{2}\mathbb{R}^{4,4}$. In particular, the
vectors $$q\wedge iq+jq\wedge kq,\quad q\wedge jq-iq\wedge
kq,\quad q\wedge kq+iq\wedge jq,\quad I_{r}p\wedge I_{s}p$$ are
contained in the orthogonal complement to $\g$.
Hence,\begin{equation}\label{(***)} R(q\wedge iq)=-R(jq\wedge
kq),\quad R(q\wedge jq)=R(iq\wedge kq),\quad R(q\wedge
kq)=-R(iq\wedge jq),\quad  R(I_{r}p\wedge I_{s}p)=0.
\end{equation}

\begin{proposition}\label{Prop1}
Any $R\in \mathcal{R}(\sp(1,1)_{\H p})$ has the following form
\begin{equation}\label{(R)}
R(I_rp\wedge
I_sq)=\Op\left(\begin{array}{cc}0&B_{rs}\\0&0\end{array}\right),\quad
R(I_rq\wedge
I_sq)=\Op\left(\begin{array}{cc}C_{rs}&D_{rs}\\0&-\bar
C_{rs}\end{array}\right),\quad R(I_{r}p\wedge I_{s}p)=0.
\end{equation} Here the numbers $C_{01},C_{02}\in \mathbb{H}$,
$d_{1},...,d_{5}\in \mathbb{R}$ are arbitrary, the numbers
$C_{rs}\in\H$, $B_{rs}, D_{rs}\in\Im\H$ are given by
\begin{align*} C_{03}&=C_{02}i-C_{01}j,\quad
C_{rs}=C_{0r}I_{s}-C_{0s}I_{r},\\\quad
B_{r0}&=\frac{1}{2}(I_1I_rC_{01}+I_2I_rC_{02}+I_3I_rC_{03}),\quad
B_{rs}=I_rC_{0s}+I_sB_{r0},
\\D_{01}&=d_{1}i+d_{2}j+d_{3}k,\quad
D_{02}=d_{2}i+d_{4}j+d_{5}k,\quad D_{03}=jD_{01}-iD_{02},\quad
D_{rs}=I_{r}D_{0s}-I_{s}D_{0r}.\end{align*}\end{proposition}

Note that $C_{rs}=-C_{sr}$, $D_{rs}=-D_{sr}$, and this is not
necessary the case for $B_{rs}$.

{\bf Proof.} Let $R\in\mathcal{R}(\sp(1,1)_{\H p})$. From
\eqref{(***)} it follows that $R(I_{r}p\wedge I_{s}p)=0.$ Define
the numbers $A_{rs}, C_{rs}\in\mathbb{H},\ B_{rs}, D_{rs}\in\Im\H$
such that
\begin{equation}\label{RRR}R(I_rp\wedge
I_sq)=\Op\left(\begin{array}{cc}A_{rs}&B_{rs}\\0&-\bar
 A_{rs}\end{array}\right),\quad
R(I_rq\wedge
I_sq)=\Op\left(\begin{array}{cc}C_{rs}&D_{rs}\\0&-\bar
C_{rs}\end{array}\right).\end{equation} From the Bianchi identity
$$R(I_{r}p,I_{s}q)I_{t}p+R(I_{s}q,I_{t}p)I_{r}p+R(I_{t}p,I_{r}p)I_{s}q=0$$
it follows $I_{t}R(I_{r}p,I_{s}q)p-I_{r}R(I_{t}p,I_{s}q)p=0$, then
$I_{t}A_{rs}p-I_{r}A_{ts}p=0$, and we have
\begin{equation}\label{(1)}
I_{t}A_{rs}-I_{r}A_{ts}=0.\end{equation}

In the same way, from
$$R(I_{r}p,I_{s}q)I_{t}q+R(I_{s}q,I_{t}q)I_{r}p+R(I_{t}q,I_{r}p)I_{s}q=0,$$
we get
\begin{equation} \label{(2)}
I_{t}B_{rs}+I_{r}C_{st}-I_{s}B_{rt}=0. \end{equation}
Also,
$$R(I_{r}q,I_{s}q)I_{t}q+R(I_{s}q,I_{t}q)I_{r}q+R(I_{t}q,I_{r}q)I_{s}q=0$$
implies that
\begin{equation}\label{(3)}
C_{rs}\bar I_{t}+C_{st}\bar I_{r}+C_{tr}\bar I_{s}=0.
\end{equation} and
\begin{equation}\label{(4)}
I_{t}D_{rs}+I_{r}D_{st}+I_{s}D_{tr}=0.
\end{equation}

From \eqref{(1)}, under the condition $t=0$ it follows that
$A_{rs}=I_{r}A_{0s}$. Substituting this  back to \eqref{(1)},
$(I_{t}I_{r}-I_{r}I_{t})A_{0s}=0$ holds. Taking $t=1$, $r=2$, we
get $A_{0s}=0$. Hence, $A_{rs}=0$ for all $r,s$.

Further,  \eqref{(3)} for $t=0$ implies
\begin{equation}\label{(3.1)}
C_{rs}=C_{0r}I_{s}-C_{0s}I_{r} \text{ for all } r,s\neq 0.
\end{equation}
From \eqref{(3)} and \eqref{(***)} for $t=1, r=2, s=3$, we have
$C_{03}=C_{02}i-C_{01}j.$ Also, from \eqref{(***)},
$C_{01}=-C_{23}, C_{02}=C_{13}, C_{03}=-C_{12}$. Similarly, using
\eqref{(4)}, we get $D_{rs}=I_{r}D_{0s}-I_{s}D_{0r}$ for all
$r,s$, from \eqref{(4)} and \eqref{(***)} $D_{03}=jD_{01}-iD_{02}$
and $D_{01}=-D_{23}, D_{02}=D_{13}, D_{03}=-D_{12}$. Note that
$D_{rs}\in \Im\H$ for any $r,s$, therefore
$D_{01}=d_{1}i+d_{2}j+d_{3}k$ and $D_{02}=d_{6}i+d_{4}j+d_{5}k$,
where $d_{1},...,d_{6}\in \mathbb{R}$. The condition $D_{03}\in
\Im\H$ implies $d_{6}=d_{2}$.

Taking $s=0$  in \eqref{(2)}, we get
$$B_{rt}=I_tB_{r0}+I_rC_{0t}.$$ Substituting this back to
\eqref{(2)}, we find $B_{r0}$. Using this, we find $B_{rs}$.
 \qed

Now let us prove that the Lie algebra $\sp(1,1)_{\H p}$ is a
Berger algebra. Consider $R$ as above, with $D_{01}=0$. Then
$R(I_{0}q\wedge
I_{1}q)=\Op\left(\begin{array}{cc}C_{01}&0\\0&-\bar
C_{01}\end{array}\right)$. Since $C_{01}$ can be arbitrary,
$\{(a,0)|a\in\mathbb{H}=\mathbb{R}\oplus\sl(1,\H)\}\subset
L(\mathcal{R}(\sp(1,1)_{\H p}))$.  On the other hand, let $R$ be
as above with $C_{rs}=0$. Since $D_{01}=d_{1}i+d_{2}j+d_{3}k$ and
$d_{1},d_{2},d_{3}$ are arbitrary, we get
$\{(0,b)|B\in\Im\H\}\subset L(\mathcal{R}(\sp(1,1)_{\H p}))$. The
inverse inclusion $L(\mathcal{R}(\sp(1,1)_{\H
p}))\subset\sp(1,1)_{\H p}$ is obvious. Thus, $\sp(1,1)_{\H
p}=L(\mathcal{R}(\sp(1,1)_{\H p}))$, i.e. $\sp(1,1)_{\H p}$ is a
Berger algebra.

Now we consider all subalgebras $\g\subset \sp(1,1)_{\H p}$, find
the spaces of curvature tensors $\mathcal{R}(\g)$  and check which
of these subalgebras are Berger algebras.

We will use the following obvious facts. Let $\g\subset \f\subset
\sp(1,1)_{\H p}$, then
\begin{equation}\label{(**)}
R\in \mathcal{R}(\g) \text{ if and only if } R\in \mathcal{R}(\f)
\text{ and } R(\wedge^2 \mathbb{R}^{4,4})\subset \g;
\end{equation}
\begin{equation}\label{(****)}
R\in \mathcal{R}(\g) \text{ if and only if } R\in \mathcal{R}(\f)
\text{ and } R|_{\g^{\perp}}=0,
\end{equation}
where $\g^{\perp}$ is the orthogonal complement to $\g$ in $\f$.

Consider a Lie subalgebra $\g\subset \sp(1,1)_{\H
p}=\mathbb{R}\oplus \sl(1,\H) \zr\Im\H$. Then we have the
following cases:

I. $\Im\H\subset\g$; \ II. $\dim(\g\cap\Im\H)=2$; \ III.
$\dim(\g\cap\Im\H)=1$;\  IV. $\g\cap\Im\H=0$.

{\bf Case I.} $\Im\H\subset\g$. In this case, $\g=\h\zr\Im\H$,
where $\h\subset\Real\oplus\sl(1,\H)$ is a subalgebra. Since there
are no 2-dimensional subalgebras of $\sl(1,\H)$, we get one of the
following subcases:

{\bf I.1.} $\h=\sl(1,\H)$. Then $\g=\sl(1,\H)\zr\Im\H$. It follows
from \eqref{(**)} that $R\in \mathcal{R}(\g)$ if and only if $R\in
\mathcal{R}(\sp(1,1)_{\H p})$ and $C_{rs}\in \Im\H$, i.e.
$C_{01}=c_{1}i+c_{2}j+c_{3}k$, $C_{02}=c_{6}i+c_{4}j+c_{5}k$,
where $c_{1},...,c_{6}\in \mathbb{R}$. The condition $C_{03}\in
\Im\H$ implies $c_{6}=c_{2}$. This shows that $\g$ is a Berger
algebra.

{\bf I.2.} $\h=\mathbb{R}(1,0)\oplus \mathbb{R}(a,0)$, where $a\in
\sl(1,\H)$.

\begin{lemma}\label{L2}  Let $(p,q)$ be an isotropic basis of $\mathbb{H}^{1,1}$ such
that $g(p,q)=1$ and $A=\Op\left(\begin{array}{cc}a&0\\0&
a\end{array}\right):\mathbb{H}^{1,1}\rightarrow\mathbb{H}^{1,1}$
is an $\mathbb{H}$-linear map, where $a=\alpha i+\beta j+\gamma
k\in\Im\H$, $\alpha,\beta,\gamma\in\mathbb{R}$ and
$\alpha^{2}+\beta^{2}+\gamma^{2}=1$. Then there exists a new
isotropic basis $(p',q')$ of $\mathbb{H}^{1,1}$ such that
$g(p',q')=1$ and  $A=\Op\left(\begin{array}{cc}i&0\\0&
i\end{array}\right)$ with respect to this basis.\end{lemma}

{\bf Proof.} Let $p'=xp$ and $q'=xq$, where
$x=x_{0}+x_{1}i+x_{2}j+x_{3}k\in\mathbb{H}$,
$x_{0},...,x_{3}\in\mathbb{R}$. Then $Ap'=xA(p)=x(\alpha i+\beta
j+\gamma k)p$. We impose the condition $Ap'=ip'$, i.e $Ap'=ixp$.
Hence, $(x_{0}+x_{1}i+x_{2}j+x_{3}k)(\alpha i+\beta j+\gamma
k)=i(x_{0}+x_{1}i+x_{2}j+x_{3}k)$. This gives the following
homogeneous system of four linear equations with respect to four
unknowns:
\begin{align*} (-\alpha+1)x_{1}-\beta x_{2}-\gamma x_{3}&=0,\quad
(\alpha-1)x_{0}+\gamma x_{2}-\beta x_{3}=0,\\
\beta x_{0}-\gamma x_{1}+(\alpha-1)x_{3}&=0,\quad \gamma
x_{0}+\beta x_{1}-(\alpha+1)x_{2}=0.\end{align*} Since the
determinant of this system equals
$(\alpha^{2}+\beta^{2}+\gamma^{2}-1)^{2}=0$, then there exists a
non-trivial solution. One can choose this solution in such a way
that $x=x_{0}+x_{1}i+x_{2}j+x_{3}k$ satisfies $x\bar x=1$. Then
$Aq'=iq'$ and $g(p',q')=g(xp,xq)=xg(p,q)\bar x=1$. \qed

By this lemma, choosing a new basis, we may assume that $a=i$.
Then $\g=(\mathbb{R}(1,0)\oplus \mathbb{R}(i,0))\zr\Im\H$. From
\eqref{(**)}, $C_{rs}\in \Real \oplus \Real i$, i.e.
$C_{01}=c_{1}+c_{2}i$ and $C_{02}=c_{3}+c_{4}i$, here
$c_{1},c_{2},c_{3},c_{4}\in \mathbb{R}$. Hence,
$C_{03}=C_{02}i-C_{01}j=c_{3}i-c_{4}-c_{1}j+c_{2}k\in
\mathbb{R}\oplus \mathbb{R}i$. This implies that $C_{01}=0$ and
$C_{02}=c_{3}+c_{4}i$, where $c_{3},c_{4}\in \mathbb{R}$ are
arbitrary. It is clear that $\g$ is a Berger algebra.

{\bf I.3.} $\h=\mathbb{R}(a,0)$, where
$a=a_{0}+a_{1}i+a_{2}j+a_{3}k\in \mathbb{R}\oplus \sl(1,\H)$,
$a_{0},...,a_{3}\in\mathbb{R}$. It follows from Lemma \ref{L2}
that there is a new basis in which $a=\alpha+\beta i$, where
$\alpha, \beta\in\mathbb{R}$ and $\alpha^{2}+\beta^{2}\neq 0$.
Thus, $\g=\mathbb{R}(\alpha+\beta i,0)\zr\Im\H$. Using
\eqref{(**)}, we have $C_{rs}\in\mathbb{R}(\alpha+\beta i)$. Let
$C_{01}=c_{1}(\alpha+\beta i)$, $C_{02}=c_{2}(\alpha+\beta i)$,
where $c_{1},c_{2}\in\mathbb{R}$. Then
$C_{03}=C_{02}i-C_{01}j=c_{2}(\alpha i-\beta)-c_{1}\alpha
j-c_{1}\beta k\in \mathbb{R}(\alpha+\beta i)$. Since
$\alpha^{2}+\beta^{2}\neq 0$, the vectors $\alpha+\beta i$ and
$-\beta+\alpha i$ are not proportional. Hence, $c_{1}=c_{2}=0$.
Thus, $C_{rs}=0$ and $\mathcal{R}(\g)=\mathcal{R}(\Im\H)$,
therefore $L(\mathcal{R}(\g))=L(\mathcal{R}(\Im\H))\neq\g$, i.e.
$\g$ is not a Berger algebra.

{\bf I.4.} $\h=0$, i.e. $\g=\Im\H$. Any $R\in\mathcal{R}(\g)$ is
given by \eqref{(R)} with $C_{rs}=B_{rs}=0$. Obviously, $\g$ is a
Berger algebra.

Note that from the two previous cases follows
\begin{lemma}\label{L1}
Any subalgebra   $\g\subset(\mathbb{R}\oplus\sl(1,\H))\zr\Im\H$
such that  $\dim\pr_{\mathbb{R}\oplus\sl(1,\H)}\g=1$  is not a
Berger algebra.
\end{lemma}

Below we will need the following lemmas:
\begin{lemma}\label{L3}  Let $(p,q)$ be an isotropic basis of
$\mathbb{H}^{1,1}$ such that $g(p,q)=1$ and in this basis
$A=\Op\left(\begin{array}{cc}0&a\\0&0\end{array}\right):\mathbb{H}^{1,1}\rightarrow\mathbb{H}^{1,1}$
is an $\mathbb{H}$-linear map, where $a=\alpha i+\beta j+\gamma
k\in\Im\H$, $\alpha^{2}+\beta^{2}+\gamma^{2}=1$. Then there exists
a new isotropic basis $(p',q')$ of $\mathbb{H}^{1,1}$ such that
$g(p',q')=1$ and
$A=\Op\left(\begin{array}{cc}0&i\\0&0\end{array}\right)$ with
respect to this basis.\end{lemma}

{\bf Proof.} As in Lemma \ref{L2}, let $p'=xp$, $q'=xq$, where
$x=x_{0}+x_{1}i+x_{2}j+x_{3}k\in\mathbb{H}$,
$x_{0},...,x_{3}\in\mathbb{R}$. Then $Aq'=xAq=xap$ and $Aq'$ is
equal to $ip'$. We get the homogeneous system as in the proof of
Lemma \ref{L2} and it has the required solution. \qed

\begin{lemma}\label{L4} Let $(p,q)$ be an isotropic basis of $\mathbb{H}^{1,1}$ such
that $g(p,q)=1$. Let
$A=\Op\left(\begin{array}{cc}0&i\\0&0\end{array}\right),
B=\Op\left(\begin{array}{cc}0&\alpha j+\beta
k\\0&0\end{array}\right):\mathbb{H}^{1,1}\rightarrow\mathbb{H}^{1,1}$
be $\mathbb{H}$-linear maps, where $\alpha,\beta\in\mathbb{R}$ and
$\alpha^{2}+\beta^{2}=1$. Then there exists a new isotropic basis
$(p',q')$ of $\mathbb{H}^{1,1}$ such that $g(p',q')=1$ and
$A=\Op\left(\begin{array}{cc}0&i\\0&0\end{array}\right)$,
$B=\Op\left(\begin{array}{cc}0&j\\0&0\end{array}\right)$
 with respect to this basis.\end{lemma}

{\bf Proof.} Let $p'=(x+yi)p$, $q'=(x+yi)q$, where
$x,y\in\mathbb{R}$. Then $Aq'=(x+yi)Aq=i(x+yi)p=ip'$. Next,
$Bq'=(x+yi)Bq=(x+yi)(\alpha j+\beta k)p$, and $Bq'$ must be equal
to $jp'$. Therefore we get the homogeneous system of two linear
equations with two unknowns:
$$(\alpha-1)x-\beta y=0,\ \beta x+(\alpha+1)y=0.$$ Since
$\alpha^{2}+\beta^{2}=1$, it follows that this system has a
non-trivial solution $(x,y)$, it can be choose in such a way that
$x^2+y^2=1$. Then
$g(p',q')=((x+yi)p,(x+yi)q)=(x+yi)g(p,q)\overline{(x+yi)}=x^{2}+y^{2}=1$.
\qed

{\bf Case II.} $\dim(\g\cap\Im\H)=2$. Using Lemmas \ref{L3} and
\ref{L4}, one can choose a new basis, such that
$\g\cap\Im\H=\mathbb{R}(0,i)\oplus\mathbb{R}(0,j)$. Then $\g$ can
be written as
$$\g=\{(a,\theta(a)k)|\ a\in
\pr_{\mathbb{R}\oplus\sl(1,\H)}\g\}\zr\mathbb{R}(0,i)\oplus\mathbb{R}(0,j),$$
where
$\theta:\pr_{\mathbb{R}\oplus\sl(1,\H)}\g\rightarrow\mathbb{R}$ is
a linear map. Let
$a=a_{0}+a_{1}i+a_{2}j+a_{3}k\in\pr_{\mathbb{R}\oplus\sl(1,\H)}\g$,
where $a_{0},...,a_{3}\in\mathbb{R}$. It holds
\begin{align*}
[(a,\theta(a)k),(0,i)]&=(0,2\Im
ia)=(0,2(a_{0}i+a_{2}k-a_{3}j))\in\g,\\
[(a,\theta(a)k),(0,j)]&=(0,2(a_{0}j-a_{1}k+a_{3}i))\in\g.\end{align*}
Hence, $a_{1}=a_{2}=0$, i.e. $a=a_{0}+a_{3}k$. This implies that
$\pr_{\mathbb{R}\oplus\sl(1,\H)}\g\subset\mathbb{R}(1,0)\oplus\mathbb{R}(k,0)$.

If $\pr_{\mathbb{R}\oplus\sl(1,\H)}\g=0$, then
$\g=\mathbb{R}(0,i)\oplus\mathbb{R}(0,j)$. Any
$R\in\mathcal{R}(\g)$ is given as in \eqref{(R)} with $C_{rs}=0$,
$B_{rs}=0$ and $D_{01}\in\mathbb{R}i\oplus\mathbb{R}j$. Also, note
that $(q\wedge kq-iq\wedge jq)\in\g^{\perp}$. On the other hand,
$(q\wedge kq+iq\wedge jq)\in(\sp(1,1)_{\mathbb{H}p})^{\perp}$ and
$\g\subset\sp(1,1)_{\mathbb{H}p}$. By \eqref{(****)}, this means
that $R(q\wedge kq)=0$, $R(iq\wedge jq)=0$, i.e. $D_{03}=0$ and
$D_{02}=-kD_{01}$. It is clear that $\g$  is a Berger algebra.

If $\pr_{\mathbb{R}\oplus\sl(1,\H)}\g=1$, then by Lemma \ref{L1},
$\g$ is not a Berger algebra.

Let $\dim\pr_{\mathbb{R}\oplus\sl(1,\H)}\g=2$, i.e.
$\pr_{\mathbb{R}\oplus\sl(1,\H)}\g=\Real(1,0)\oplus\Real(k,0)$. In
the case $\theta=0$ we have
$\g=(\mathbb{R}(1,0)\oplus\mathbb{R}(k,0))\zr(\mathbb{R}(0,i)\oplus\mathbb{R}(0,j))$.
Let us check that $\g$ is a Berger algebra. Let $R\in\R(\g)$. As
it is observed above, $(q\wedge kq-iq\wedge jq)\in\g^{\perp}$ and
$(q\wedge kq+iq\wedge jq)\in\g^{\perp}$. By \eqref{(****)}, this
means that $R(q\wedge kq)=0$, $R(iq\wedge jq)=0$, i.e. $D_{03}=0$
and $D_{02}=-kD_{01}$. Moreover, from \eqref{(**)},
$D_{01}\in\mathbb{R}i\oplus\mathbb{R}j$. For the same reason,
$C_{03}=0$ and $C_{02}=C_{01}k$. Also, from \eqref{(**)},
$C_{01}\in \mathbb{R}\oplus\mathbb{R}k$. Thus $\g$ is a Berger
algebra.

Further, suppose that $\theta\neq 0$. Let
$\ker\theta=\mathbb{R}(\alpha+\beta k)$, where $\alpha,
\beta\in\mathbb{R}$ and $\alpha^{2}+\beta^{2}\neq 0$. Then
$$\g=(\mathbb{R}(\alpha+\beta k,0)\oplus\mathbb{R}(-\beta+\alpha
k,\gamma k))\zr(\mathbb{R}(0,i)\oplus\mathbb{R}(0,j)),$$ where
$\gamma=\theta (-\beta+\alpha k).$ Next, $$[(\alpha+\beta
k,0),(-\beta+\alpha k,\gamma k)]=(0,2\alpha\gamma k)\in\g.$$
Hence, $\alpha=0$. Then we may assume that $\beta=1$.
  With respect to the new basis $p'=p$, $q'=q+\frac{\gamma}{2}kp$, we get
$\g=(\mathbb{R}(1,0)\oplus\mathbb{R}(k,0))\zr(\mathbb{R}(0,i)\oplus\mathbb{R}(0,j))$.
This algebra has been already considered above.

{\bf Case III.} $\dim(\g\cap\Im\H)=1$. By Lemma \ref{L3}, in the
new basis $\g\cap\Im\H=\mathbb{R}(0,i)$. Then
$$\g=\{(a,\theta(a))|\
a\in\pr_{\mathbb{R}\oplus\sl(1,\H)}\g\}\zr\mathbb{R}(0,i),$$ where
$\theta:\pr_{\mathbb{R}\oplus\sl(1,\H)}\g\rightarrow\mathbb{R}(0,j)\oplus\mathbb{R}(0,k)$
is a linear map. Let
$a=a_{0}+a_{1}i+a_{2}j+a_{3}k\in\pr_{\mathbb{R}\oplus\sl(1,\H)}\g$,
where $a_{0},...,a_{3}\in\mathbb{R}$. Since $\g$ is a Lie algebra,
we get $[(a,\theta(a)),(0,i)]=(0,2(a_{0}i+a_{2}k-a_{3}j))\in\g$.
Hence, $a_{2}=a_{3}=0$, i.e.
$\pr_{\mathbb{R}\oplus\sl(1,\H)}\g\subset\mathbb{R}(1,0)\oplus\mathbb{R}(i,0)$.
One can easily see that if $\pr_{\mathbb{R}\oplus\sl(1,\H)}\g=0$,
i.e $\g=\mathbb{R}(0,i)$, then $\g$ is not a Berger algebra. Also,
if $\dim\pr_{\mathbb{R}\oplus\sl(1,\H)}\g=1$, then by Lemma
\ref{L1}, $\g$ is not a Berger algebra.

Now let $\dim\pr_{\mathbb{R}\oplus\sl(1,\H)}\g=2$. Then
$\pr_{\mathbb{R}\oplus\sl(1,\H)}\g=\mathbb{R}(1,0)\oplus\mathbb{R}(i,0)$.
Writing $\theta(a)=\theta_{1}(a)j+\theta_{2}(a)k$, where
$a\in\mathbb{R}\oplus\mathbb{R}i$, we obtain
\begin{multline*}[(1,\theta_{1}(1)j+\theta_{2}(1)k),(i,\theta_{1}(i)j+\theta_{2}(i)k)]=(0,2\Im(\theta_{1}(i)j+\theta_{2}(i)k))-(0,2\Im
(\theta_{1}(1)j+\theta_{2}(1)k)i)\\=(0,2(\theta_{2}(i)+\theta_{1}(1))k+2(\theta_{1}(i)
-\theta_{2}(1))j)\in\g.\end{multline*} Hence,
$\theta_{2}(1)=\theta_{1}(i)$ and  $\theta_{2}(i)=-\theta_{1}(1)$.
There exist $\alpha, \beta\in\mathbb{R}$ such that
$\theta(1)=-\alpha j+\beta k$, $\theta(i)=\beta j+\alpha k$.
Consequently, $\g=\mathbb{R}(1,-\alpha j+\beta
k)\oplus\mathbb{R}(i,\beta j+\alpha k)\zr\mathbb{R}(0,i)$. Using
Mathematica 4.0, we find that $\mathcal{R}(\g)=0$, i.e. $\g$ is
not a Berger algebra.

{\bf Case IV.}  $\dim(\g\cap\Im\H)=0$. Then $\g=\{(a,\theta(a))|\
a\in\pr_{\mathbb{R}\oplus\sl(1,\H)}\g\}$, where
$\theta:\pr_{\mathbb{R}\oplus\sl(1,\H)}\g\rightarrow\Im\H$ is a
linear map. Put
$\theta(a)=\theta_{1}(a)i+\theta_{2}(a)j+\theta_{3}(a)k$. As
follows from Lemma \ref{L1}, if
$\dim\pr_{\mathbb{R}\oplus\sl(1,\H)}\g=1$, then $\g$ is not a
Berger algebra. Now consider the case when
$\dim\pr_{\mathbb{R}\oplus\sl(1,\H)}\g\geq 2$, i.e.
$\pr_{\mathbb{R}\oplus\sl(1,\H)}\g\cap\sl(1,\H)\neq 0$. Let
$a\in\pr_{\mathbb{R}\oplus\sl(1,\H)}\g\cap\sl(1,\H)$ and
$\theta(a)\neq 0$. By Lemma \ref{L3}, there is a new basis such
that $\theta(a)=i$. Then we have $A=(a,i)$. Choosing again a new
basis $p'=p$, $q'=q-\frac{a}{\|a\|}ip$, one can get $A=(a,0)$.
Further, from Lemma \ref{L2} it follows that $A=(i,0)$, i.e.
$\theta(i)=0$. Since $\g$ is a Lie algebra,
\begin{align*}
[(i,0),(a,\theta(a))]&=(2a_{3}j-2a_{2}k,2\theta_{3}(a)j-2\theta_{2}(a)k)\in\g,\\
[(i,0),(a_{3}j-a_{2}k,\theta_{3}(a)j-\theta_{2}(a)k)]&=
(-2a_{2}j-2a_{3}k,-2\theta_{2}(a)j-2\theta_{3}(a)k)\in\g,\end{align*}
here $a=a_{0}+a_{2}j+a_{3}k$, i.e. $(a_{0}+a_{2}j+a_{3}k,
\theta_{1}(a)i+\theta_{2}(a)j+\theta_{3}(a)k)\in\g$ and
$(a_{2}j+a_{3}k,\theta_{2}(a)j+\theta_{3}(a)k)\in\g$. Hence,
$(a_{0},\theta_{1}(a)i)\in\g$. Since $(a_{0},
\theta_{1}(a_{0})i+\theta_{2}(a_{0})j+\theta_{3}(a_{0})k)\in\g$,
we get $\theta_{2}(a_{0})=\theta_{3}(a_{0})=0$ and
$\theta_{1}(a_{2}j+a_{3}k)=0$, i.e.
$\theta_{1}(a)=\theta_{1}(a_{0})$.

There are two cases:

{\bf 1.} $\pr_{\mathbb{R}\oplus\sl(1,\H)}\g\subset\sl(1,\H)$, i.e.
$a_{0}=0$ for any $a\in\pr_{\mathbb{R}\oplus\sl(1,\H)}\g$. We
already know that $(0,i)\in\g$ and all other elements from $\g$
have the form
$(a_{2}j+a_{3}k,\theta_{2}(a_{2}j+a_{3}k)j+\theta_{3}(a_{2}j+a_{3}k)k)$.
Since $\g$ is a Lie algebra, then
$$[(i,0),(a_{2}j+a_{3}k,\theta_{2}(a_{2}j+a_{3}k)j+\theta_{3}(a_{2}j+a_{3}k)k)]=(2a_{2}k-2a_{3}j,2\theta_{2}(a_{2}j+a_{3}k)k-2\theta_{3}(a_{2}j+a_{3}k)j)\in\g.$$
On the other hand,
$\theta(a_{3}j-a_{2}k)=\theta_{2}(a_{3}j-a_{2}k)j+\theta_{3}(a_{3}j-a_{2}k)k$.
Hence, $\theta_{2}(a_{3}j-a_{2}k)=\theta_{3}(a_{2}j+a_{3}k)$ and
$\theta_{3}(a_{3}j-a_{2}k)=-\theta_{2}(a_{2}j+a_{3}k)$. Since
$a_{2},a_{3}\in\mathbb{R}$ are arbitrary,  these equations are
equivalent to $\theta_{2}(k)=-\theta_{3}(j)$ and
$\theta_{3}(k)=\theta_{2}(j)$. Let $\theta_{2}(j)=\alpha$ and
$\theta_{2}(k)=\beta$. Thus,
$\g=\mathbb{R}(i,0)\oplus\mathbb{R}(j,\alpha j-\beta
k)\oplus\mathbb{R}(k,\beta j+\alpha k)$. We check with Mathematica
4.0 that $\dim\mathcal{R}(\g)=0$, i.e. $\g$ is not a Berger
algebra.

{\bf 2.} There is $a\in\pr_{\mathbb{R}\oplus\sl(1,\H)}\g$ such
that $a_{0}\neq 0$, i.e. $\g\nsubseteq\sl(1,\H)\zr\Im\H$. From
above, $(1,\theta_{1}(1)i)\in\g$ and
$\g=\mathbb{R}(1,\theta_{1}(1)i)\oplus\mathbb{R}(i,0)\oplus\mathbb{R}(j,\alpha
j-\beta k)\oplus\mathbb{R}(k,\beta j+\alpha k)$. We have
$[(1,\theta_{1}(1)i),(j,\alpha j-\beta  k)]=(0,2\alpha j-2\beta
k-2\theta_{1}(1)k)\in\g$ and $\g\cap\Im\H=0$. Therefore,
$\alpha=0$ and $\theta_{1}(1)=-\beta$. Consequently,
$$\g=\mathbb{R}(1,-\beta i)\oplus\Real(i,0)\oplus\Real(j,-\beta k)\oplus\Real(k,\beta
j).$$ Choosing the basis $p'=p,q'=q-\frac{1}{2}\beta ip$, we get
$\g=\{(a,0)|a\in\H\}$. Let $R\in\R(\g)$, then $D_{rs}=0$ and
$B_{rs}=0$. Consequently, $C_{rs}=0$ and $R=0$, i.e. $\g$ is not a
Berger algebra.
 The theorem is proved.  \qed

\section{Proof of Theorem \ref{Th3}}
Let $\g\subset\sp(1)\oplus\sp(1,1)_{\mathbb{H} p}$ and
$\g\not\subset\sp(1,1)_{\mathbb{H} p}$. From \eqref{RR0W}
 it follows that \begin{equation}\label{Rsp1+}
\R(\sp(1)\oplus\sp(1,1))=\Real
R_0\oplus\R(\sp(1,1)),\end{equation} where $R_0$ is given by
\eqref{R0} with $X,Y\in\Real^{4,4}$. Consider the subalgebra
\begin{equation}\h_0=\sp(1)\oplus\{(a,0)|a\in\Real\oplus\sl(1,\H)\}=\left.
\left\{\left(\begin{array}{cc}A&0\\0&-A^t\end{array}\right)\right|A\in\Real
E_4\oplus\so(4)\right\}\subset \sp(1)\oplus\sp(1,1)_{\mathbb{H}
p}.\end{equation} It is not hard to check that the element $R_1$
defined by
\begin{equation}\label{R1} R_1(I_sp,I_rq)=-I_rp\wedge I_sq+I_sp\wedge I_r q-\delta_{sr}\left(\begin{array}{cc}E_4&0\\0&-E_4\end{array}\right), \quad
R_1(I_sp,I_rp)=R_1(I_sq,I_rq)=0\end{equation} belongs to the space
$\R(\h_0)$. From this, \eqref{Rsp1+} and the equality
$\R(\{(a,0)|a\in\Real\oplus\sl(1,\H)\})=0$ it follows that
$\R(\h_0)=\Real R_1$. This and \eqref{Rsp1+} show that
\begin{equation}\label{Rsp1++} \R(\sp(1)\oplus\sp(1,1)_{\mathbb{H} p})=\Real R_1\oplus\R(\sp(1,1)_{\mathbb{H} p}).\end{equation}

Suppose that $\g\subset\sp(1)\oplus\sp(1,1)_{\mathbb{H} p}$ is a
Berger algebra,  $\g\not\subset\sp(1,1)_{\mathbb{H} p}$ and
$\g\neq\h_0$. Then there exists $R\in\R(\g)$ such that $R=\nu
R_1+R'$, $\nu\neq 0$ and $R'\in\R(\sp(1,1)_{\H p})$. From above we
get that $\pr_{\h_0} R(I_rp,I_sq)=R_1(I_rp,I_sq)$, hence
$\pr_{\h_0}\g=\h_0$, here the projection is taken with respect to
the decomposition $$\sp(1)\oplus\sp(1,1)_{\mathbb{H} p}=\h_0\zr
\Im\H.$$ Hence for some $b\in \Im \H$ the element $\xi=(1,b)$
belongs to $\g$. If $b\neq 0$, then choosing the basis $p,q'$,
where $q'=q-\frac{1}{2}bp$, we get that $\xi=(1,0)\in\g$. Let
$\xi_1=A+(a,b)\in\g$, where $A\in\sp(1)$,
$a\in\Real\oplus\sl(1,\H)$ and $b\in \Im \H$. Then,
$[\xi,\xi_1]=(0,2b)\in\g$. This shows that $\g=\h_0\zr L$, where
$L\subset\Im\H$. The Lie brackets of elements from $\h_0$ and $L$
are given by the representation of $\Real\oplus\sl(1,\H)$ on $\Im
\H$. Since this representation is irreducible,
$\g=\sp(1)\oplus\sp(1,1)_{\H p}$. Thus if
$\g\subset\sp(1)\oplus\sp(1,1)_{\mathbb{H} p}$ is a Berger algebra
and  $\g\not\subset\sp(1,1)_{\mathbb{H} p}$, then either
$\g=\sp(1)\oplus\sp(1,1)_{\H p}$, or $\g=\h_0$. Obviously, the
both algebras are Berger algebras. \qed

\section{The spaces $\R(\g)$}\label{secR}

In the proofs of the previous theorems we have found the spaces
$\R(\g)$ for all Berger subalgebras
$\g\subset\sp(1)\oplus\sp(1,1)_{\H p}$. Here we summaries these
results. For the algebra $\g_1=\sp(1,1)_{\H p}$ the space
$\R(\g_1)$ is found in Proposition \ref{Prop1}. For the algebras
$\g_2,...,\g_6$ the spaces $\R(\g)$ can be found as in Proposition
\ref{Prop1} with the following additional conditions:
\begin{description}\item[$\g_2$:] $C_{01}=c_1i+c_2j+c_3k$,
$C_{02}=c_2i+c_4j+c_5k$, where $c_1,...,c_5\in\Real$ are
arbitrary;
\item[$\g_3$:]
$C_{01}=0$, $C_{02}=c_0+c_1i$, where $c_0,c_1\in\Real$ are
arbitrary;
\item[$\g_4$:]
$C_{01}=C_{02}=0$;
\item[$\g_5$:]
$C_{01}=c_0+c_3k$, $C_{02}=c_0k-c_3$, $d_3=d_5=0$, $d_4=-d_1$,
where $c_0,c_3,d_1,d_2\in\Real$ are arbitrary;
\item[$\g_6$:]
$C_{01}=C_{02}=0$, $d_3=d_5=0$, $d_4=-d_1$, where
$d_1,d_2\in\Real$ are arbitrary.
\end{description}
Finally, $\R(\sp(1)\oplus\sp(1,1)_{\H p})=\Real
R_1\oplus\R(\sp(1,1)_{\H p})$ and $\R(\h_0)=\Real R_1$, where
$R_1$ is given by \eqref{R1}.

%%%%%%%%%%%%%%%%%%%%%%%%%%%%%%%%%%%%%%

\section{Pseudo-quaternionic-K\"ahlerian symmetric spaces of signature $(4,4)$}\label{secSym} Simply connected pseudo-Riemannian symmetric spaces with reductive
holonomy algebras are classified by Berger in \cite{Ber57}. In
\cite{AC05} it is shown that pseudo-quaternionic-K\"ahlerian
symmetric spaces of non-zero scalar curvature have reductive
holonomy algebras and all these spaces were listed. In
\cite{AC01,KO07} indecomposable simply connected
pseudo-hyper-K\"ahlerian symmetric spaces of signature $(4,4n)$
are classified. Here we use the results of this paper to give new
proof to the classification of simply connected
pseudo-quaternionic-K\"ahlerian symmetric spaces of signature
$(4,4)$ with non-irreducible holonomy algebras, in particular, we
find the holonomy algebras and the curvature tensors of these
spaces.

First we give some preliminaries. Recall that any simply
connected pseudo-Riemannian symmetric space $(M,g)$ of signature
$(r,s)$ is uniquely defined by a triple $(\h,\sigma,\eta)$ ({\it a
symmetric triple}), where $\h$ is a Lie algebra, $\sigma$ is an
involutive automorphism of $\h$ and $\eta$ is a non-degenerate
$\ad_{\h_+}$-invariant  symmetric bilinear form of signature
$(r,s)$ on the vector space  $\h_-$.  Here $\h_\pm$ are the
eigenspaces of $\sigma$ corresponding to the eigenvalues $\pm 1$.
One may assume that $[\h_-,\h_-]=\h_+$. The vector space $\h_-$
can be identified with the tangent space to $M$ at a fixed point
$o\in M$. Then the form $\eta$ is identified with the form $g_o$.
The curvature tensor $R_o$ of the manifold $(M,g)$ is given by
$$R_o(x,y)z=-[[x,y],z],$$
where $x,y,z\in\h_-=T_oM$. For the holonomy algebra $\g$ of the
manifold $(M,g)$ at the point $o\in M$ we have
$$\g=R_o(T_oM,T_oM)=[\h_-,\h_-]=\h_+.$$ Since $\nabla R=0$, $\g$ annihilates $R_o\in\R(\g)$, i.e.
\begin{equation}\label{AR=0}A\cdot R_o=0,\quad (A\cdot R_o)(x,y)= [A,R_o(x\wedge y)]-R_o(Ax\wedge y)
-R_o(x\wedge Ay),\end{equation} where $A\in\g$ and $x,y\in
\Real^{r,s}$. For a subalgebra $\g\subset\so(r,s)$ denote by
$\R_0(\g)$ the subset of $\R(\g)$ consisting of curvature tensors
$R$ satisfying \eqref{AR=0}. If $(M,g)$ is a symmetric space with
the curvature tensor  $R$ and holonomy algebra $\g$, then
$R_o\in\R_0(\g)$, and $R_o(\Real^{r,s},\Real^{r,s})=\g$, where
$o\in M$. Conversely, let $\g\subset\so(r,s)$ be a subalgebra and
$R\in\R_0(\g)$ with $R(\Real^{r,s},\Real^{r,s})=\g$. Consider the
Lie algebra $\h=\g+\Real^{r,s}$ with the Lie brackets
$$[x,y]=-R(x,y),\quad [A,x]=Ax,\quad [A,B]=A\circ B-B\circ A,$$
where $x,y\in\Real^{r,s}$ and $A,B\in\g$. Define the involutive
automorphism $\sigma$ of $\h$ by $\sigma(A+x)=A-x,$ where $A\in
\g$, $x\in \Real^{r,s}$. We have $\h_+=\g$ and $\h_-=\Real^{r,s}$.
Obviously, the pseudo-Euclidean metric $\eta$ on $\Real^{r,s}$ is
$\ad_{\h_+}$-invariant. We get the symmetric triple
$(\h,\sigma,\eta)$ that defines a simply connected
pseudo-Riemannian symmetric space of signature $(r,s)$. Thus any
pseudo-Riemannian symmetric space $(M,g)$ of signature $(r,s)$
defines a pair $(\g,R)$ ({\it a symmetric pair}), where
$\g\subset\so(r,s)$ is a subalgebra, $R\in\R_0(\g)$ and
$\g=R(\Real^{r,s},\Real^{r,s})$ holds. Conversely, any such pair
defines a simply connected pseudo-Riemannian symmetric space.

{\it An isomorphism }of  symmetric triples
$f:(\h_1,\sigma_1,\eta_1)\to (\h_2,\sigma_2,\eta_2)$ is a Lie
algebra isomorphism $f:\h_1\to\h_2$ such that
$f\circ\sigma_1=\sigma_2\circ f$ and the induced map
$f:\h_{1-}\to\h_{2-}$ is an isometry of the pseudo-Euclidean
spaces $(\h_{1-},\eta_1)$ and $(\h_{2-},\eta_2)$. Isomorphic
symmetric triples define isometric pseudo-Riemannian symmetric
spaces. {\it An isomorphism }of  symmetric pairs $f:(\g_1,R_1)\to
(\g_2,R_2)$ consists of an isometry of $\Real^{r,s}$ that defines
the equivalence of the representations $\g_1,\g_2\subset\so(r,s)$
and sends $R_1$ to $R_2$. Obviously, two symmetric pairs are
isomorphic if and only if they define isometric symmetric triples.
For a positive real number $c\in\Real$, the symmetric pairs
$(\g,R)$ and $(\g,cR)$ define diffeomorphic simply connected
symmetric spaces and the metrics of these spaces differ by the
factor $c$. Note that the symmetric pairs $(\g,R)$ and $(\g,-R)$
can define non-isometric symmetric spaces, e.g. in the Riemannian
case if $\g$ is irreducible. In the case of neutral signature
$(r,r)$, we may pass from  $(\g,R)$ to $(\g,-R)$  multiplying the
metric by $-1$.

\begin{theorem} Let $(M,g)$ be a non-flat simply connected pseudo-hyper-K\"ahlerian symmetric
space of signature $(4,4)$ and $\g\subset\sp(1,1)$ its holonomy
algebra. Then $\g$ preserves an isotropic 4-dimensional
$I_1,I_2,I_3$-invariant subspace of $\Real^{4,4}$, and $(M,g)$ is
given by exactly one of the following symmetric pairs $(\g,R)$:
\begin{itemize} \item[(1)] \begin{align*}
\g&=\left.\left\{\Op\left(\begin{array}{cc}0&b\\0&0\end{array}\right)\right|b\in\Im\H\right\},\\
R(q,I_1q)&=-R(I_2q,I_3q)=\Op\left(\begin{array}{cc}0&
i\\0&0\end{array}\right),\qquad
R(q,I_2q)=R(I_1q,I_3q)=\Op\left(\begin{array}{cc}0&\lambda
j\\0&0\end{array}\right),\\
R(q,I_3q)&=-R(I_1q,I_2q)=\Op\left(\begin{array}{cc}0&-(1+\lambda)k\\0&0\end{array}\right),\quad
R(I_rp,\cdot)=0,\end{align*} where
$\lambda\in\left[-\frac{1}{2},0\right)\cup(0,+\infty)$ is a
parameter.
\item[(2)]\begin{align*} \g&=\left.\left\{\Op\left(\begin{array}{cc}0&b\\0&0\end{array}\right)\right|b\in\Real i\oplus\Real
j\right\},\qquad
R(q,I_1q)=-R(I_2q,I_3q)=\Op\left(\begin{array}{cc}0&
i\\0&0\end{array}\right),\\
R(q,I_2q)&=R(I_1q,I_3q)=\Op\left(\begin{array}{cc}0&-
j\\0&0\end{array}\right),\quad R(q,I_3q)=R(I_1q,I_2q)=
R(I_rp,\cdot)=0.\end{align*}
\end{itemize}\end{theorem}

{\bf Proof.}  It is well known that there are no Ricci-flat
symmetric spaces with irreducible holonomy algebras, hence $\g$ is
not irreducible. First for each possible holonomy algebra
$\g\subset\sp(1,1)_{\H p}$ we find the elements of $R\in\R_0(\g)$
such that $R(\Real^{4,4},\Real^{4,4})=\g$.

Suppose that $\g$ contains $A=(1,0)$. Let $R\in\R(\g)$ given by
\eqref{RRR}. Assume that $\g$ annihilates $R$. Then $$0=(A\cdot
R)(I_r q,I_sq)=(0,2D_{rs})+2R(I_r q,I_sq)=2(C_{rs},2D_{rs}),$$
hence, $R=0$ and it does not span $\g$.

Further, consider $\g=\{(a,b)|\, a,b\in\Im\H\}$. Let $A=(a,0)$,
$a\in\Im\H$ and $R\in\mathcal{R}(\g)$ is given as in Proposition
1. Then $$(A\cdot
R)(I_{r}q,I_{s}q)=\Op\left(\begin{array}{cc}C_{rs}a-aC_{rs}&2\Im
D_{rs}a\\0&C_{rs}a-aC_{rs}\end{array}\right)-R(I_{r}aq,I_{s}q)-R(I_{r}q,I_{s}aq).$$
 Since any $R\in\mathcal{R}(\g)$ is given by
$C_{01},C_{02}$ and $D_{01},D_{02}$ it is enough to consider only
the cases $r=0$, $s=1$ and $r=0$, $s=2$. Let $a=i$. We obtain
\begin{equation}\label{(SU1)} C_{01}i-iC_{01}=0, \quad
C_{02}i-iC_{02}=-2C_{03}.
\end{equation} For the last equation we also used \eqref{(***)}.
Similarly, taking $a=j$, $r=0,s=1$ and  $r=0,s=2$, we get
\begin{equation}\label{(SU2)} C_{01}j-jC_{01}=2C_{03}, \quad
C_{02}j-jC_{02}=0.
\end{equation} From the first equations of
\eqref{(SU1)} and \eqref{(SU2)} it follows that $C_{01}=c_{1}i$,
$C_{02}=c_{2}j$, where $c_{1},c_{2}\in\mathbb{R}$. Recall that
$C_{03}=C_{02}i-C_{01}j$, hence $C_{03}=-(c_2+c_{1})k$. From this
and the second equations of \eqref{(SU1)} and \eqref{(SU2)} we get
$c_{1}=c_{2}=0$. Consequently, $C_{rs}=0$ and $B_{rs}=0$ for all
$r,s$. Hence $R$ does not span $\g$.

On the other hand, $\R_0(\g)=\R(\g)$ for the Lie algebras
$\g=\{(0,b)|b\in\Im\H\}$ and $\g=\{(0,b)|b\in\Real i\oplus\Real
j\}$. Thus only these two Lie algebras  may appear as the holonomy
algebras of  symmetric spaces.

Let $\g=\{(0,b)|b\in\Im\H\}$. We will check which $R\in\R(\g)$
define different symmetric spaces. Recall that any $R\in\R(\g)$ is
uniquely given by the numbers $d_1,...,d_5\in\Real$ in the
following way: $$R(I_rp,\cdot)=0,\quad R(I_rq,I_sq)=(0,D_{rs}),$$
where $D_{03}=jD_{01}-iD_{02}$, $D_{rs}=I_{r}D_{0s}-I_{s}D_{0r}$,
$D_{01}=d_{1}i+d_{2}j+d_{3}k$, $D_{02}=d_{2}i+d_{4}j+d_{5}k$. This
shows that $R$ is defined by the values $R(q,yq)$, $y\in\Im\H$.
Define the $\Real$-linear map $\varphi:\Im \H\to\Im\H$ by the
equality
$R(q,yq)=\Op\left(\begin{array}{cc}0&\varphi(y)\\0&0\end{array}\right)$.
For its matrix with respect to the basis $i,j,k$ we get
$\Mat\varphi=\left(\begin{array}{ccc}d_1& d_2&d_3\\d_2&
d_4&d_5\\d_3& d_5&-d_1-d_4\end{array}\right)$.  We see that the
matrix $\Mat\varphi$ is symmetric and trace-free. The condition
$R(\Real^{4,4},\Real^{4,4})=\g$ is equivalent to
$\det\Mat\varphi\neq 0$.

Let $(\g,R)$ be a symmetric pair as above. Consider the new basis
$(p'=xp,q'=xq)$, where $x\in\H$ and $x\bar x=1$. We get
$R(q',yq')=\Op\left(\begin{array}{cc}0&\varphi_1(y)\\0&0\end{array}\right)$,
where $y\in\Im\H$, $\varphi_1:\Im \H\to\Im\H$ is a linear map, and
the matrix is considered with respect to the basis $(p',q')$. Then
$R(q',yq')=\Op\left(\begin{array}{cc}0&\bar{x}\varphi_1(y)x\\0&0\end{array}\right)$
 with respect to the basis $(p,q)$.
 On the other hand,
$$R(q',yq')=R(xq,yxq)=R(q,\bar{x}yxq)=\Op\left(\begin{array}{cc}0&\varphi(\bar{x}yx)\\0&0\end{array}\right),$$
where the matrix is considered with respect to the basis $(p,q)$.
We conclude that $\bar{x}\varphi_1(y)x=\varphi(\bar{x}yx)$. Hence,
$\varphi_1(y)=x\varphi(\bar{x}yx)\bar{x}$. Recall that the group
of the unit quaternionic numbers  acts on $\Im\H$ as the group of
special orthogonal transformations by the rule $x\cdot
y=A_xy=\bar{x}yx$. Thus, $$\Mat\varphi_1=(\Mat
A_x)^{-1}\cdot\Mat\varphi\cdot \Mat A_x,$$ where the matrices are
considered with respect to the basis $(i,j,k)$ of $\Im\H$.
Consequently, there exists $x$ such that $\Mat\varphi_1$ is a
diagonal matrix. Let $(\mu,\lambda,\upsilon)$ be its diagonal
elements. Note that the triple $(\mu,\lambda,\upsilon)$ is defined
up to a permutation. By the above, $\mu\lambda\upsilon\neq 0$ and
$\mu+\lambda+\upsilon=0$. We may assume that $\mu>0$. Since
$(\g,R)$ and  $(\g,cR)$ ($c>0$) define the same symmetric space,
we may take $\mu=1$. Obviously, two symmetric pairs are equivalent
if and only if  they define the same unordered triples
$(1,\lambda,\upsilon)$. Due to the symmetry between
$\lambda,\upsilon$, we get that any symmetric pair $(\g,R)$ is
equivalent to exactly one symmetric pair given by the ordered
triple $(1,\lambda,-1-\lambda)$, where
$\lambda\in\left[-\frac{1}{2},0\right)\cup(0,+\infty)$. This
proves the statement for the first algebra.

Let $\g=\{(0,b)|b\in\Real i\oplus\Real j\}$. In this case any
$R\in\R(\g)$ is defined   by the value
$R(q,I_1q)=\Op\left(\begin{array}{cc}0&
D_{01}\\0&0\end{array}\right)$, where $D_{01}=d_1i+d_2j$
$d_1,d_2\in\Real$. Such $R$ defines a linear map $\varphi:\Real
i\oplus\Real j\to\Real i\oplus\Real j$. It is symmetric and
trace-free. As above, we may consider the new basis
$(p'=xp,q'=xq)$, where $x\in\H$ and $x\bar x=1$. We must assume
that $\g$ has the same matrix form with respect to the new basis,
then either $x=x_0+kx_3$, or $x=x_1i+x_2j$. We may find a basis
such that $\varphi$ is given by a diagonal matrix with the
diagonal elements $(\mu,-\mu)$, the rescalling gives $\mu=1$. Thus
all symmetric pairs with this holonomy algebra define isometric
simply connected symmetric spaces. \qed

\begin{theorem} Let $(M,g)$ be a  simply connected pseudo-quaternionic-K\"ahlerian symmetric
space of signature $(4,4)$ with non-irreducible holonomy algebra
$\g\subset\sp(1)\oplus\sp(1,1)$ and non-zero scalar curvature.
Then
$$\g=\sp(1)\oplus\left.\left\{\Op\left(\begin{array}{cc}a&0\\0&-\bar
a\end{array}\right)\right|a\in\H\right\}$$ and
 $(M,g)$ is given by $(\g,R_1)$,
 where $R_1$ is defined by \eqref{R1}.
\end{theorem}

Note that the obtained symmetric space is isometric to ${\rm
SL}(2,\H)/{\rm S}({\rm GL}(1,\H)\times {\rm GL}(1,\H))$. This
follows from the classification of \cite{AC05} and from the fact
that ${\rm SL}(n+1,\H)/{\rm S}({\rm GL}(1,\H)\times {\rm
GL}(n,\H))$ are the only non-flat simply connected
pseudo-quaternionic-K\"ahlerian symmetric spaces with
non-irreducible holonomy algebras \cite{AC05,Ber57}.

{\bf Proof.} As in the previous theorem it is easy to prove that
the first algebra of Theorem \ref{Th3} satisfies $\R_0(\g)=\Real
R_1$, where $R_1$ is defined by \eqref{R1}. The image of this
tensor does not span $\g$. For the second $\g$ it holds
$\R_0(\g)=\Real R_1$. In this case, $R_1$ do span $\g$.  \qed

\end{document}